\documentclass[12pt]{article}
\usepackage{amsmath}
\usepackage{amscd}
\usepackage{graphicx}
\usepackage{amssymb,amsfonts,amsthm,amscd,mathrsfs}
\usepackage{verbatim}
\usepackage[T2A]{fontenc}
\usepackage[utf8]{inputenc}
\usepackage[russian]{babel}
\usepackage{bm}
\usepackage{authblk}

\renewcommand{\b}{\beta}

\newcommand{\e}{\varepsilon}

\newcommand{\E}{\mathsf{E}}

\begin{document}

\title{On Symmetrized Pearson's Type Test for Normality of Autoregression:  Power under Local Alternatives
}

\author{M.~V.~Boldin    
\footnote{Moscow State Lomonosov Univ., Dept. of Mech. and Math., Moscow, Russia\\
e-mail: boldin$_{-}$m@hotmail.com}}

\date{ }
\maketitle

\textbf{Abstract}

We consider a stationary linear AR($p$) model with observations subject to
gross errors  (outliers). The autoregression parameters as well as the distribution function (d.f.) $G$ of innovations are unknown. The distribution of outliers $\Pi$ is unknown and arbitrary, their intensity is
 $\gamma n^{-1/2}$ with an unknown $\gamma$, $n$ is the sample size.
 We test the hypothesis for normality of innovations $$\bm H_{\Phi}\colon G \in \{\Phi(x/\theta),\,\theta>0\},$$ $\Phi(x)$ is the d.f. $\bm N(0,1)$.  
Our test is the special symmetrized Pearson's type test. We find the power of this test under local alternatives 
$$\bm H_{1n}(\rho)\colon G(x)=A_n(x):=(1-\rho n^{-1/2})\Phi(x/\theta_0)+\rho n^{-1/2}H(x),
$$ $\rho\geq 0,\,\theta_0$ is the unknown (under $\bm H_{\Phi}$) variance of innovations.
 First of all we estimate the autoregression parameters  
  and then using the residuals from
the estimated autoregression we construct a kind of empirical distribution function
(r.e.d.f.), which is a counterpart of the (inaccessible)
e.d.f. of the autoregression innovations. After this we construct the symmetrized variant r.e.d.f. Our test statistic is the functional from symmetrized r.e.d.f. We obtain a stochastic expansion
of this  symmetrized r.e.d.f.  under $\bm H_{1n}(\rho)$ , which enables us to  investigate our test.
We establish qualitative robustness of this test in terms of uniform equicontinuity
 of the limiting power  (as functions of $\gamma,\rho$ and $\Pi$) with respect to $\gamma$ in a neighborhood of $\gamma=0$.

{\bf Key words:} autoregression, outliers, residuals, empirical distribution function,
Pearson's chi-square test, robustness, estimators, normality, local alternatives.

{\bf 2010 Mathematics Subject Classification:} Primary 62G10; secondary 62M10, 62G30, 62G35.

\section{
Введение }
В этой работе мы рассматриваем стационарную AR($p$) модель с ненулевым средним
\begin{equation}
v_t = \b_1 v_{t-1} + \dots + \b_p v_{t-p} +\nu + {\e}_t, \quad  t \in \mathbb{Z} .
\end{equation}

В (1.1) $\{{\e}_t\}$ -- независимые одинаково распределенные случайные величины (н.о.р.сл.в.) с неизвестной функцией распределения (ф.р.) $G(x)$;
 $\E \e_1 = 0$,  $0<\E {\e}_1^2 < \infty$; $\bm{\b} = (\b_1, \dots, \b_p)^T \in  \mathbb{R}^p $-- вектор неизвестных параметров, таких что корни соответствующего (1.1) характеристического уравнения по модулю меньше единицы;  $\nu$ -- неизвестное среднее,  $\nu \in \mathbb R^1$.\\
 Эти условия дальше всегда предполагаются выполненными и особо не оговариваются.\\
 Мы рассматриваем модель (1.1) с выбросами в наблюдениях. А именно, предполагается, что наблюдаются величины
$$
y_t=v_t+z^{\gamma_n}_t {\xi}_t,\quad t  = 1-p, \dots, n, \eqno(1.2)
$$
где $v_{1-p}, \dots, v_n$ -- выборка из стационарного решения $\{v_t\}$ уравнения  (1.1); $\{z^{\gamma_n}_t\}$ н.о.р.сл.в., имеющие распределение Бернулли, т.е. принимающие значения 1 и 0, причем вероятность единицы
 $\gamma_n$,
 $$
 \gamma_n = \min(1, \frac{\gamma}{\sqrt{n}}) ,\quad  \gamma \ge 0 \text{\; неизвестно.}
 $$
Кроме того, $\{ \xi_t \}$ -- н.о.р.сл.в. с произвольным и неизвестным распределением $\Pi$.
Переменные $\{ \xi_t \}$ интерпретируются как выбросы (засорения), $\gamma_n$ уровень засорения. Для $\gamma = 0$ мы получаем модель (1.1) без засорений.

Модель (1.2) -- локальный вариант хорошо известной модели засорения данных во временных рядах, см. \cite{MartYoh86}.\\
В этой работе мы рассматриваем гипотезу о нормальности инноваций 
$$
\bm{H}_{\Phi}\colon G(x) \in \{\Phi(x/\theta),\; \theta > 0\},
$$
где $\Phi(x)$ -- стандартная нормальная ф.р. Напомним, что гипотеза $H_{\Phi}$ эквивалентна нормальности самой стационарной последовательности 
 $\{v_t\}$.\\
   Нормальность инноваций обеспечивает оптимальность процедур наименьших квадратов оценивания и проверки гипотез в авторегрессии, см., например, \cite{And.},\cite{Broc.Dav.}. Поэтому проверка $\bm H_{\Phi}$ -- содержательная задача. 
 В \cite{Bold.Arx.2} была построена специальная симметризованная статистика типа хи-квадрат Пирсона для $\bm H_{\Phi}$ в схеме (1.1) --(1.2) с неизвестным средним и найдено ее распределение при гипотезе. Необходимый следующий шаг -- найти мощность построенного теста при локальных альтернативах и исследовать устойчивость мощности к выбросам. Этот шаг мы и делаем в настоящей работе.\\ 
 Обозначим неизвестную дисперсию инноваций при $\bm H_{\Phi}$ через $\theta_0$.  
  В качестве альтернативы к $\bm H_{\Phi}$ берется предположение о том, что $\{\e_t\}$  в (1.1) -- н.о.р.сл.в. с ф.р. в виде смеси
 $$ 
 G(x)=A_n(x):=(1-\rho_n)\Phi(x/\theta_0)+\rho_nH(x),\eqno(1.3)
 $$
 где
 $$
  H(x)-\text{ф.р.}, \rho_n=\min(1,\frac{\rho}{\sqrt{n}}), \rho\geq 0.
 $$
 Предположение (1.3) будем понимать как локальную альтернативу к $\bm H_{\Phi}$ и обозначать $\bm{H}_{1n}(\rho)$. Разумеется, $\bm H_{\Phi}$ и $\bm H_{1n}(0)$ совпадают.\\ Наша цель -- найти предельное распределение симметризованной статистики из \cite{Bold.Arx.2} при $\bm{H}_{1n}(\rho)$, найти мощность теста и исследовать ее устойчивость к выбросам.\\ 
   Отметим, что недавно в \cite{Bold.Petr.},\cite{Bold.2020} были получены результаты о проверке  
   $\bm H_{\Phi}$ в схеме (1.1)--(1.2)  при $\nu=0$. Результаты получены и при гипотезе и при локальных альтернативах. Но нулевое среднее -- существенное ограничение.  В нашей ситуации при неизвестном и, быть может, ненулевом $\nu$, использовать статистики типа хи квадрат из \cite{Bold.Petr.}, \cite{Bold.2020} не удается, т.к. их предельное распределение зависит от оценок $\bm{\b}$ и $\nu$ даже при $\gamma=0$.\\
 Мы преодолеваем эту трудность,  
 строя симметризованную статистику Пирсон. Разбиение, которым она определяется, специальное и симметрично относительно нуля. Симметризованная статистика Пирсона  является функционалом от симметризованной остаточной эмпирической функции распределения (о.э.ф.р.).  Мы находим асимптотические распределения симметризованной о.э.ф.р. и нашей статистики Пирсона при  $\bm{H}_{1n}(\rho)$. Асимптотические распределения симметризованных о.э.ф.р. и статистики Пирсона при $\bm H_{\Phi}$ и  $\gamma=0$ свободны.  
 Кроме того, симметризованный тест Пирсона оказывается качественно робастным.
 
 Все определения и результаты (основные, это Теоремы 2.1--2.2) представлены в Разделе 2.
 
 \section{Основные результаты}

\subsection{Стохастическое разложение о.э.ф. р. }

Относительно альтернативы $\bm H_{1n}(\rho)$ из (1.3) нам потребуется сделать некоторые предположения.

{\bf Условие (i).} Случайная величина с функцией распределения $H(x)$ имеет нулевое средние и конечную дисперсию.

{\bf Условие (ii).} Функции распределения $H(x)$ дважды дифференцируема с ограниченной второй производной.

{\bf Условие (iii).} Функции распределения $H(x)$ симметрична относительно нуля.

Перепишем уравнение (1.1) в удобном для дальнейшего рассмотрения виде. Для этого определим константу $\mu$ соотношеним
$$
\nu=(1-\b_1-\ldots-\b_p)\mu,
$$
тогда
$$
v_t -\mu= \b_1 (v_{t-1}-\mu) + \dots + \b_p( v_{t-p} -\mu)+ {\e}_t, \quad  t \in \mathbb{Z.}
$$

Если положить $u_t:=v_t-\mu,$ то
$$
v_t=\mu+u_t,\quad u_t=\b_1 u_{t-1} + \dots + \b_pu_{t-p}+ {\e}_t, \quad  t \in \mathbb{Z.}\eqno(2.1	)
$$
Последовательность $\{u_t\}$ в (2.1) -- авторегрессионная последовательность с нулевым средним и конечной дсперсией.\\
Построим по наблюдениям $\{y_t\}$ из (1.2) оценки ненаблюдаемых $\{\e_t\}$.\\ 
Далее $\Gamma,\, R$ --любые конечные неотрицательные числа. Пусть $\hat\mu_n$ будет любая последовательность, для которой при $\bm H_{1n}(\rho)$ последовательнсть
$$
n^{1/2}(\hat\mu_n-\mu)=O_P(1),\quad n\to\infty, \text{ равномерно по}\,\,\rho\leq R,\,\,\gamma \leq\Gamma.\eqno(2.2)
$$
 Положим
$$
\hat u_t=y_t-\hat\mu_n,\quad t=1,\ldots,n.
$$
Пусть $\hat{\bm\b}_n= (\hat \b_{1n}, \dots, \hat \b_{pn})^T$ будет любая последовательность,
для которой при $\bm H_{1n}(\rho)$ последовательность\\
 $$
 n^{1/2}(\hat{\bm\b}_n-\bm\b)=O_P(1),\quad n\to \infty, \text{ равномерно по}\,\,\rho\leq R,\,\,\gamma\leq \Gamma.\eqno(2.3)
 $$
Примеров подходящих оценок $\hat\mu_n,\,\hat{\bm\b}_n$ много, широкий класс составляют, например, М-оценки, построенные по $\{y_t\}$ и $\{\hat u_t\}$ аналогично тому, как они строятся по незасоренным данным, см., например,   \cite{Koul87}. Результаты для М-оценок при альтернативе $\bm H_{1n}(\rho)$ в схеме (1.1)--(1.2) вполне аналогичны изложенным в Разделе 2.4 \cite{Bold.Petr.}. \\
Положим
$$
\hat \e_t =  \hat {u}_t -  \hat \b_{1n}\hat { u}_{t-1} - \dots - \hat \b_{pn}\hat{ u}_{t-p},\quad t = 1, \dots, n.\eqno(2.4)
$$
Величины $\{\hat \e_t\}$ называются остатками. 

Ведем остаточную эмпирическую функцию распределения ( о.э.ф.р.)
$$
\hat G_n(x) = n^{-1} \sum_{t=1}^n I(\hat \e_t \le x),\quad  x\in \mathbb{R}^1.
$$
 Здесь и в дальнейшем $I(\cdot)$ обозначает индикатор события.\\
Функция $\hat G_n(x)$ -- аналог гипотетической э.ф.р.
$$
G_n(x) = n^{-1} \sum_{t=1}^n I(\e_t \le x)
$$
ненаблюдаемых величин $\e_{1}, \dots, \e_n$.\\ 
Следствие 2.1 в \cite{Bold.Arx.3} прямо влечет следующую Теорему 2.1. 
\newtheorem{Th}{Теорема}[section]
\begin{Th}
Предположим, что верна альтернатива $\bm H_{1n}(\rho)$. Пусть выполнены Условия (i) -- (ii), и пусть $\varphi(x)$ будет стандартная гауссовская плотность. Пусть $\Gamma \ge 0,\,R\ge 0$ будут любые конечные числа. 
Тогда для любого $\delta > 0$
$$
\sup_{\rho\leq R,\gamma \le \Gamma} \Prob(|n^{1/2}[\hat G_n(x) -G_n(x)]- \frac{1}{\theta_0}\varphi(\frac{x}{\theta_0})\hat\delta_n n^{1/2}(\hat\mu_n-\mu)-  \gamma\Delta_0(x, \Pi) | > \delta) \to 0,
\quad n \to \infty.
$$

Здесь сдвиг
$$
\Delta_0(x, \Pi) = \sum_{j=0}^p [\E G_0(x + \b_j \xi_1) - G_0(x)],\, \b_0 =-1,\,G_0(x)=\Phi(\frac{x}{\theta_0}),
$$
а $\hat\delta_n=1-\hat \b_{1n}-\ldots-\hat \b_{pn}$.
\end{Th}

При гипотезе $\bm H_{\Phi}$  ф.р. $G(x)=\Phi(\frac{x}{\theta_0})$ и, значит, симметрична относительно нуля. Поэтому будем брать оценкой $G(x)$ симметризованную оценку
$$
\hat S_n(x):=\frac{\hat G_n(x)+1-\hat G_n(-x)}{2}.\eqno(2.5)
$$
Положим
$$
\Delta_S(x, \Pi):=\frac{\Delta_0(x, \Pi)-\Delta_0(-x, \Pi)}{2}.\eqno(2.6)
$$
Пусть
$$
 S_n(x):=\frac{ G_n(x)+1- G_n(-x)}{2}.\eqno(2.7)
$$

Теорема 2.1 и четность $\varphi(x)$  прямо влекут
\newtheorem{Corollary}{Следствие
}[section]
\begin{Corollary}
При условиях Теоремы 2.1 
$$
\sup_{\rho\leq R,\gamma \le \Gamma} \Prob(|n^{1/2}[\hat S_n(x) -S_n(x)]-  \gamma\Delta_S(x, \Pi) | > \delta) \to 0,
\quad n \to \infty.
$$
\end{Corollary}

\subsection{Симметризованный тест типа хи-квадрат Пирсона для $\bm H_{\Phi}$}

В этом разделе мы построим тест типа хи-квадрат Пирсона для проверки 
$$
\bm H_{\Phi}\colon  G(x) \in \{\Phi(x/\theta),\;\theta>0\},
$$
где $\Phi(x)$ -- стандартная нормальная ф.р. Как отмечалось во Введении, такая гипотеза эквивалентна гауссовости стационарного решения (1.1). Наша цель -- найти распределение тестовой статистики при альтернативе $\bm H_{1n}(\rho)$ 
 из (1.3).

Напомним, $\theta_0$ истинное и неизвестное значение $\theta$ при $\bm H_{\Phi}$, тогда при $\bm H_{\Phi}$ $G(x)=\Phi(x/\theta_0)$.
Для полуинтервалов 
$$
B^+_j=(x_{j-1},x_j],\quad j=1,\ldots,m,\; m>2,\; 0=x_0<x_1<\ldots<x_m=\infty,
$$
пусть 
$$
p^+_j(\theta)=\Phi(x_{j}/\theta)-\Phi(x_{j-1}/\theta)>0.
$$
При $\bm H_{\Phi}$ $\Prob(\e_1 \in B^+_j)=p^+_j(\theta_0)$. Если ввести еще симметричные полуинтервалы

 $ B^-_j=(-x_{j},-x_{j-1}]$, то при $\bm H_{\Phi}$
$$
\Prob(\e_1 \in B^+_j \cup B^-_j)=2p^+_j(\theta_0):=p_j(\theta_0),\,\sum_{j=1}^m p_j(\theta_0)=1.
$$
Пусть $\hat \nu_j^+$ обозначает число остатков среди $\{\hat \e_t,\,t=1,\ldots,n\}$, попавших в $B^+_j$, а 
$\hat \nu_j^-$ обозначает число остатков, попавших в $B^-_j$. 
Тогда
$$
\hat \nu_j^+=n[\hat G_n(x_j)-\hat G_n(x_{j-1})],\quad 
\hat \nu_j^-=n[\hat G_n(-x_{j-1})-
\hat G_n(-x_{j})].\eqno(2.8)
$$
Пусть
$$
\hat\nu_j=\hat \nu_j^+ +\hat \nu_j^-.
$$
Оценкой $\theta_0$ мы возьмем $n^{1/2}$-состоятельное решение уравнения
$$
\sum_{j=1}^m \frac{\hat{\nu}_j}{p_j(\theta)}p_j^{'}(\theta)=0, 
$$
в котором штрих означает производную относительно  $\theta$. Такое уравнение -- аналог обычного уравнения модифицированного метода хи-квадрат, в котором ненаблюдаемые частоты $\nu_1,\ldots,\nu_m$ инноваций, попавших в $B_1,\ldots,B_m,\,B_j:=B_j^+\cup B_j^-$, заменены их оценками $\hat{\nu}_1,\ldots,\hat{\nu}_m$.

Интересующая нас тестовая статистика типа хи-квадрат для $\bm H_{\Phi}$ имеет вид
$$
\hat{\chi}^2_{n} = \sum_{j=1}^m \frac{(\hat{\nu}_j - n p_j(\hat{\theta}_n))^2}{n p_j(\hat{\theta}_n)}.\eqno(2.9)
$$
Статистика $\hat{\chi}^2_{n} $ является функционалом от $\hat S_n(x)$. Действительно, в силу (2.5) и (2.8)
$$
n[\hat S_n(x_j)-\hat S_n(x_{j-1})]=\frac{1}{2}\{n[\hat G_n(x_j)-\hat G_n(x_{j-1})]+n[\hat G_n(-x_{j-1})-\hat G_n(-x_{j})]\}=
$$
$$
(\hat\nu^+_j +\hat\nu^-_j)/2=\hat\nu_j/2.\eqno(2.10)
$$
Будем называть $\hat{\chi}^2_{n}$ из (2.9) симметризованной статистикой Пирсона.\\
Чтобы описать асимптотические свойства $\hat{\chi}^2_{n}$ нам понадобятся некоторые обозначения. Далее $ |\cdot|$ означает Евклидову норму вектора или матрицы.

Введем диагональную матрицу
$$
\bm P_0=\mbox{diag}\{p_1(\theta_0),\ldots,p_m(\theta_0)\}
$$
и векторы
$$
\bm p_0=(p_1(\theta_0),\ldots,p_m(\theta_0))^T,\quad
\bm p_0'=(p_1'(\theta_0),\ldots,p_m'(\theta_0))^T,
$$
$$
\bm b_0=\bm P_0^{-1/2}\bm p'_0, \quad \bm \alpha_0=\bm b_0/|\bm b_0|,
$$
$$
\bm p_A=(p_1^A,\ldots,p_m^A)^T \text{с компонентами}\, p_j^A=2(A_n(x_j)-A_n(x_{j-1})),
$$

$$
\bm p_H=(p_1^H,\ldots,p_m^H)^T,\quad p_j^H:=2(H(x_j)-H(x_{j-1})).
$$
Нам понадобится еще вектор
$$
\bm \delta(\Pi)=(\delta_1(\Pi),\ldots,\delta_m(\Pi))^T\quad \text{с компонентами}\,\, \delta_j(\Pi)=2[\Delta_{S}(x_{j},\Pi)-\Delta_{S}(x_{j-1},\Pi)],
$$
сдвиг $\Delta_{S}(x,\Pi)$ определен в (2.6).

Из определения $S_n(x)$ и $\nu_j$ аналогично (2.10) следует:
$$
n[ S_n(x_j)- S_n(x_{j-1})]=\nu_j/2.\eqno(2.11)
$$
В силу соотношений (2.10) --(2.11), Следствия 2.1 и определения вектора $\delta(\Pi)$ имеем при $\bm H_{1n}(\rho)$:
$$
n^{1/2}(\frac{\hat\nu_j}{n}-
p_j(\theta_0))-n^{1/2}(\frac{\nu_j}{n}-p_j(\theta_0))=
n^{1/2}(\frac{\hat\nu_j}{n}-\frac{\nu_j}{n})=
$$
$$
2n^{1/2}\{[\hat S_n(x_j)-\hat S_n(x_{j-1}]-[S_n(x_j)-S_n(x_{j-1})]\}=\gamma \delta_j(\Pi)+o_P(1),\quad n\to \infty.\eqno(2.12)
$$
Здесь $o_P(1)$ обозначает величину, сходящуюся к нулю по вероятности равномерно по $\rho\leq R,\,\gamma \leq \Gamma$.\\
Введем вектора 
$$
\hat {\bm \nu}=(\hat\nu_1,\ldots,\hat\nu_m)^T,\quad \bm \nu=(\nu_1,\ldots,\nu_m)^T.
$$
В силу (2.12) при $\bm H_{1n}(\rho)$
$$
n^{1/2}(\frac{\hat{\bm \nu}}{n}-\bm p_0)=
n^{1/2}(\frac{\bm \nu}{n}-\bm p_A)+\rho(\bm p_H-\bm p_0)+\gamma \bm\delta(\Pi)+o_P(1),\quad n\to \infty.\eqno(2.13)
$$
Напомним еще известный факт относительно слабой сходимости вектора $\bm \nu$ (см. \cite{Bold.2020},  соотношение (3.10)): при $\bm H_{1n}(\rho)$ равномерно по $\rho \leq R$
$$
n^{1/2}(\frac{\bm \nu}{n}-\bm p_A) \to \bm N(\bm 0,\,\, \bm P_0-\bm p_0\bm p_0^T),\quad n\to \infty.\eqno(2.14)
$$
Соотношений (2.13) --(2.14) достаточно, чтобы следуя схеме доказательства Теоремы 2.2 в \cite{Bold.2020} доказать наше  основное утверждение -- Теорему 2.2. В ней $ F_{m-2}(x, \hat\lambda^2(\rho,\gamma, \Pi))$ означает ф.р. нецентрального хи-квадрата с $m-2$ степенями свободы и параметром нецентральности $\hat\lambda^2(\rho,\gamma, \Pi)$, а $\bm E_m$ означает единичную матрицу порядка  $m$.
\begin{Th}
Пусть альтернатива  $\bm H_{1n}(\rho)$ верна. Пусть выполнены Условия (i)--(iii). Тогда для любых конечных  неотрицательных $\Gamma,\,R$
$$
\sup_{x\in\mathbb{ R}^1, \rho\leq R,\gamma \le \Gamma} |\Prob(\hat{\chi}^2_{n} \le x) - F_{m-2}(x, \hat\lambda^2(\rho,\gamma, \Pi))| \to 0,\quad n \to \infty.
$$
Параметр нецентральности равен
$$
\hat\lambda^2(\rho,\gamma, \Pi) =
|(\bm E_m-\bm \alpha_0 \bm  \alpha_0^T) \bm{ P}^{-1/2}_0[\rho(\bm p_H-\bm p_0)+\gamma \bm{\delta}(\Pi)]|^2.
$$
\end{Th}
В силу Теоремы 2.2 при $H_{\Phi}$, т.е. при $\rho=0$, и $\gamma=0$ предельное распределение нашей тестовой статистики $\hat{\chi}^2_{n}$ будет обычное (центральное) распределение хи-квадрат с $m-2$ степенями свободы.

Для заданного $0<\alpha<1$ ($\,\alpha$ --асимптотический уровень в схеме без засорений), мы будем отвергать $\bm H_{\Phi}$, когда
$$
\hat{\chi}^2_{n} > \chi_{m-2}(1-\alpha), \eqno(2.15)
$$
$\chi_{m-2}(1-\alpha)$ -- квантиль уровня $1-\alpha$ распределения хи-квадрат с $m-2$ степенями свободы. 
Мощность такого теста равна
$$
W_n(\rho,\gamma,\Pi)=P(\hat{\chi}^2_{n} > \chi_{m-2}(1-\alpha)).
$$
В силу Теоремы 2.2 эта мощность сходится при $n\to\infty$ равномерно по $\rho\leq R,\,\,\gamma\leq \Gamma$
к асимптотической мощности
$$
W(\rho,\gamma,\Pi)= 1 - F_{m-2}(\chi_{m-2}(1-\alpha), \hat\lambda^2(\rho,\gamma, \Pi)),\quad
 W(0,0,\Pi)= \alpha.
 $$
 Используя неравенство
 $$
 |F_k(x,\lambda_1^2)-F_k(x,\lambda_2^2)|\leq 2 \sup_{x \in \mathbb{ R}^1}\varphi(x)|\lambda_1-\lambda_2|,
 $$
 $\varphi(x)$ стандартная гауссовская плотность, и определение $\hat\lambda^2(\rho,\gamma, \Pi))$, получаем:
 $$
\sup_{\Pi,\rho\geq 0}|W(\rho,\gamma,\Pi)-W(\rho,0,\Pi)|\leq \gamma \sqrt{2/\pi}|(\bm E_m-\bm \alpha_0 \bm  \alpha_0^T) \bm{ P}^{-1/2}_0|\sup_{\pi}| \bm {\delta}( \Pi)| \to 0,\quad \gamma \to 0.\eqno(2.16)
$$
Соотношение (2.16) означает асимптотическую качественную робастность теста (2.15). Такая робастность означает, что  $\bm H_{\Phi}$ для малых  $\gamma$ можно проверять примерно с асимптотическим уровнем $\alpha$, и тест будет иметь примерно такую же асимптотическую мощность при $\bm H_{1n}(\rho)$ c $\rho>0$, как в схеме без засорений. И все это  независимо от распределения засорений $\Pi$.

\end{document}